\pdfoutput=1

\documentclass[ECP]{ejpecp} 

\usepackage[utf8]{inputenc}

\usepackage{enumitem}
\usepackage{cleveref}

\SHORTTITLE{Characterisation of gradient flows on finite state Markov chains}

\TITLE{Characterisation of gradient flows
  on~finite~state~Markov~chains\thanks{Support: UK Engineering and Physical
    Sciences Research Council (EPSRC) grant EP/H023348/1 for the University of
    Cambridge Centre for Doctoral Training, the Cambridge Centre for
    Analysis.}}

\AUTHORS{%
  Helge~Dietert\footnote{University of Cambridge, UK.
    \EMAIL{H.G.W.Dietert@maths.cam.ac.uk}}}


\AMSSUBJ{60J27}

\SUBMITTED{May 12, 2014}
\ACCEPTED{March 1, 2015}

\ARXIVID{1405.2552v2}

\VOLUME{20}
\YEAR{2015}
\PAPERNUM{29}
\DOI{v20-3521}

\ABSTRACT{
  In his 2011 work, Maas has shown that the law of any time-reversible
  continuous-time Markov chain with finite state space evolves like a
  gradient flow of the relative entropy with respect to its stationary
  distribution. In this work we show the converse to the above by
  showing that if the relative law of a Markov chain with finite state
  space evolves like a gradient flow of the relative entropy
  functional, it must be time-reversible.
  When we allow general functionals in place of the relative entropy,
  we show that the law of a Markov chain evolves as gradient flow if
  and only if the generator of the Markov chain is real
  diagonalisable. Finally, we discuss what aspects of the functional
  are uniquely determined by the Markov chain.
}

\newcommand{\dd}{\mathrm{d}}

\newcommand{\spX}{\ensuremath{\mathcal{X}}} 
\newcommand{\dP}{\ensuremath{\mathcal{P}}} 
\newcommand{\R}{\ensuremath{\mathbb{R}}} 
\newcommand{\fF}{\ensuremath{\mathfrak{F}}} 
\newcommand{\relH}{\ensuremath{\mathfrak{H}}} 

\begin{document}



\section{Introduction}
The seminal paper of Jordan, Kinderlehrer, and Otto
\cite{jordan1998variational} identified Markov processes in the
continuous setting as gradient flows of the entropy using the
Wasserstein distance. This understanding lead to many new results (see
Villani~\cite{villani2009optimal} for an overview).  More recently,
Maas \cite{maas2011gradient} considered Markov chains with finite
state space and showed that, in this case, the Wasserstein distance
does not allow this identification. Instead, assuming a
time-reversible Markov chain, he was able to construct a different
metric that allows this identification. Different constructions have
been given in \cite{chow2012fokker,mielke2011gradient} and the setting
used is also described in \cite{erbar2014gradient}.

The construction of the metric is involved and uses time-reversibility
at several places. This further motivates our study of the converse of
these statements. For this we will first introduce the setting used
and define the gradient flow.

We consider continuous-time irreducible Markov chain with finite state
space $\spX = \{0,1,\dots ,d\}$. We denote its generator by
$Q \in \R^{\spX \times \spX}$ where for $i\not = j$ the entry $Q_{ij}$
is the transition rate from state $i$ to state $j$ and
$Q_{ii} = -\sum_{j\not = i} Q_{ij}$.

Given the initial probability distribution $\mu$ of the Markov chain,
the probability distribution after time $t$ will be given by
$\mu e^{tQ}$. Note that the transition matrix $e^{tQ}$ acts on the
right on the row vector $\mu$ and the evolution of the law is captured
by the Markov semigroup $e^{tQ}$. Since the Markov chain is
irreducible, there exists a unique stationary distribution $\pi$ to
which $\mu e^{tQ}$ will converge as $t \to \infty$.

For the definition of a gradient flow, let $\dP$ be the space of
probability distributions on $\spX$ with positive mass for any
state. Then $\dP$ can be naturally understood as the $d$-dimensional
sub-manifold
$\{v\in\R^\spX : \sum_{i\in\spX} v_i = 1 \text{ and } v_i > 0\;
\forall i\in\spX\}$
of $\R^\spX$. Under this identification, the tangent space at any point
is
\begin{equation*}
  T = \{ v\in\R^\spX : \sum_{i\in\spX} v_i = 0 \}.
\end{equation*}

Given a functional $\fF : \dP \mapsto \R$ and a Riemannian metric $g$
on $\dP$, the gradient flow $\rho : \R^+\mapsto \dP$ is determined by
the differential equation
\begin{equation*}
  g|_{\rho(t)}(\dot{\rho}(t),v) = - \dd \fF|_{\rho(t)}(v) \quad \forall v \in T,
  t \in \R^+.
\end{equation*}
That is $\dot\rho$ equals $-\dd \fF$ under the identification of the
tangent space and the cotangent space through the metric $g$.

We say that the gradient flow under $g$ of $\fF$ equals the flow
associated to the semigroup $e^{tQ}$ if for all $\mu \in \dP$ the
trajectory $t \mapsto \mu e^{tQ}$ equals the gradient flow $\rho(t)$
with $\rho(0) = \mu$.

From the continuous setting a natural functional is the relative
entropy $\relH$ with respect to the stationary state $\pi$, which is
defined by
\begin{equation*}
  \relH(\rho) = - \sum_{i\in\spX} \rho_i \log \frac{\pi_i}{\rho_i}.
\end{equation*}
Throughout this work the relative entropy is understood with respect
to the stationary distribution of the considered Markov chain.

The result by Maas \cite{maas2011gradient} now is: \emph{There exists
  a metric $g$ on $\dP$ such that the gradient flow under $g$ of the
  relative entropy $\relH$ equals the flow associated to the semigroup
  $e^{tQ}$.}

First we show that the metric is not unique:
\begin{theorem}
  \label{thm:pert}
  For $d \ge 2$, consider a continuous-time irreducible Markov chain
  with finite state space $\spX = \{0,1,\dots ,d\}$, generator $Q$,
  and stationary distribution $\pi$. Let $\fF : \dP \mapsto \R$ be
  a differentiable functional and $g$ a Riemannian metric on $\dP$.
  If the flow associated to the semigroup $e^{tQ}$ equals the gradient
  flow of $\fF$ under the metric $g$, then for $\rho\in\dP$ with
  $\rho\not = \pi$, there exists another metric $\tilde{g}$ on $\dP$
  such that $\tilde{g} \not = g$ at $\rho$ and such that the gradient
  flow of $\fF$ under the metric $\tilde{g}$ still equals the flow
  associated to the semigroup $e^{tQ}$.
\end{theorem}
As converse of the construction we show:
\begin{theorem}
  \label{thm:pi}
  Consider a continuous-time irreducible Markov chain with finite
  state space $\spX = \{0,1,\dots ,d\}$, generator $Q$, and
  stationary distribution $\pi$.
  If the flow associated to the semigroup $e^{tQ}$ equals the gradient
  flow of $\fF\in C^2$ under a Riemannian metric $g\in C^1$ on $\dP$,
  then, with $g|_\pi$ as metric at $\pi$,
  \begin{enumerate}[label=(\alph*), ref=(\alph*)]
  \item \label{item:g-unique} $g|_\pi$ is uniquely determined by $Q$ and $\fF$,
  \item \label{item:q-dependence} $Q$ can be computed from $g|_\pi$ and $\fF$,
  \item \label{item:q-real-diag} $Q$ is real diagonalisable,
  \item \label{item:q-time-reversible} $Q$ is time-reversible if $\fF$
    is the relative entropy $\relH$.
  \end{enumerate}

  Conversely, if $Q$ is real diagonalisable, then there exists a
  Riemannian metric $g$ and a smooth functional $\fF$ on $\dP$ such
  that the gradient flow of $\fF$ under the metric $g$ equals the flow
  associated to the semigroup $e^{tQ}$.
\end{theorem}
This result shows that the assumption of time-reversibility in the
construction of the metric by Maas is necessary and cannot be
relaxed. Moreover, the results for the generator $Q$ come from the
differentiability around the equilibrium distribution $\pi$, so that
the theorem holds as long as there is a neighbourhood of $\pi$ in
which the Markov semigroup $e^{tQ}$ equals the gradient flow.
\begin{remark}
  \Cref{thm:pert,thm:pi} are obtained by analysing the Riemannian
  structure of the gradient flow and thus can be formulated for
  general gradient flows on finite-dimensional manifolds. For this, part
  \ref{item:q-time-reversible} of \Cref{thm:pi} can be formulated with
  a weighted $\ell^2$-norm, i.e. $Q$ is symmetric with respect to this
  norm, if $\fF$ is the squared distance to $\pi$ under this
  $\ell^2$-norm. In fact, relating the results to $Q$ acting on the
  bigger space $\R^{\spX}$ makes the proofs slightly longer.
\end{remark}

Combining this theorem with Maas' result gives our main theorem.
\begin{theorem}[Characterisation of Markov chains]
  \label{thm:char}
  Consider a continuous-time irreducible Markov chain with finite
  state space and generator $Q$.
  \begin{itemize}
  \item The Markov chain is time-reversible if and only if there
    exists a metric $g$ such that the flow associated to the semigroup
    $e^{tQ}$ equals the gradient flow under $g\in C^1$ of the relative
    entropy with respect to the stationary distribution.
  \item The Markov chain has a real diagonalisable generator $Q$ if
    and only if there exists a metric $g\in C^1$ and a functional
    $\fF\in C^2$ such that the flow associated to the semigroup
    $e^{tQ}$ is the gradient flow of $\fF$ under $g$.
  \end{itemize}
\end{theorem}

The characterisation of real diagonalisable generators $Q$ shows for
example that Markov chains that also have an oscillatory behaviour
cannot be described by gradient flows.

Finally, we remark that the relative entropy $\relH$ depends on the
generator $Q$ only through the equilibrium distribution
$\pi$. Moreover, any functional $\fF$ that allows to construct a
gradient flow for all time-reversible Markov chains must have a
similar Taylor expansion around $\pi$. More precisely:
\begin{theorem}
  \label{thm:relH}
  Fix a finite state space $\spX$, a distribution
  $\pi\in\dP$ and a functional $\fF:\dP\mapsto \R$ in $C^2$.
  Suppose that, for every generator $Q$ defining an irreducible
  time-reversible Markov chain with state space $\spX$ and stationary
  distribution $\pi$, there exists a Riemannian metric $g\in C^1$ on
  $\dP$ such that the gradient flow of $\fF$ under the metric $g$
  equals the flow associated to the semigroup $e^{tQ}$. Then there
  exists a positive constant $\alpha$ such that
  \begin{equation*}
    \dd \fF|_{\pi} = \dd \relH|_{\pi} = 0
    \text{ and }
    \dd^2 \fF|_{\pi} = \alpha\, \dd^2 \relH|_{\pi},
  \end{equation*}
  where $\relH$ is the relative entropy with respect to $\pi$.
\end{theorem}

Here we use the notation $\dd \fF|_{\pi}$ to denote the first
derivative of $\fF$ at the point $\pi$, which we understand as linear
map from $T$ to $\R$. With $\dd^2 \fF|_{\pi}$ we denote the second
derivative at the point $\pi$, which is a linear map from $T \times T$
to $\R$.

The assumption on the functional $\fF$ is not empty, because Maas'
result states that the relative entropy $\relH$ with respect to $\pi$
is a functional satisfying the assumption. Moreover, it cannot be
strengthened to uniqueness. For this, another functional is the
quadratic form $\fF$ defined by $\fF|_{\pi} = \dd \fF|_{\pi} = 0$ and
$\dd^2 \fF|_{\pi} = \alpha \dd^2 \relH|_{\pi}$ for some $\alpha >
0$.
This satisfies the assumption, because for any generator $Q$ the
constant metric $g$ given by the value $g|_{\pi}$ in part
\ref{item:g-unique} of \Cref{thm:pi} indeed defines a Riemannian
metric with the required identification.

\begin{remark}
  In \cite{maas2011gradient}, Maas considered continuous-time Markov
  chains obtained from an irreducible discrete-time Markov chain with
  transition matrix $K$ by choosing the jump times according to a
  Poisson process. The resulting generator is $Q=K-I$ and, by analogy
  with continuous-time diffusion processes, the semigroup $e^{tQ}$ is
  also called a heat flow.

  By time-rescaling, this is no restriction to the class of Markov
  chains for the study of gradient flows because, for any generator
  $Q$, we can find some $\alpha > 0$ such that $K=I+\alpha Q$ is
  non-negative along the diagonal and $K$ defines a transition
  matrix. Now given a functional $\fF$, if we can find a metric $g$
  such that the flow associated to the semigroup $e^{t(\alpha Q)}$
  equals the gradient flow of $\fF$ under $g$, then the flow
  associated to the semigroup $e^{tQ}$ equals the gradient flow of
  $\fF$ under the rescaled metric $g/\alpha$.
\end{remark}

\section{Characterisation of Markov chains}
Recall that for an irreducible Markov chain the evolution $\mu
e^{tQ}$ converges to the unique stationary distribution $\pi$ for any
$\mu\in\dP$. Hence $\mu Q$ vanishes if and only if
$\mu=\pi$. With this observation, we can construct a perturbation of
the metric in order to prove \Cref{thm:pert}.

\begin{proof}[Proof of \Cref{thm:pert}]
  Let $e_1$ be the vector field given by $\mu Q$ at $\mu\in\dP$. For a
  small enough neighbourhood $V \subset \dP$ of $\rho$, we can find
  smooth vector fields $e_2,\dots ,e_d$ such that $e_1,\dots ,e_d$ is
  a basis of $T$ at every $\mu\in V$. Let $\eta$ be a smooth cutoff
  functional with compact support, vanishing outside $V$, and
  satisfying $\eta(\rho)\not = 0$. Then define another metric
  $\tilde{g}$ by
  \begin{equation*}
    \tilde{g}|_{\mu}(e_i,e_j) =
    \begin{cases}
      g(e_i,e_j)|_{\mu} + \eta(\mu) a & \text{if $i=j=2$,}\\
      g(e_i,e_j)|_{\mu} & \text{otherwise,}
    \end{cases}
  \end{equation*}
  for $\mu\in V$ and a constant $a\in\R$. Outside of $V$, define
  $\tilde{g}=g$.

  If $a$ is small enough, $\tilde{g}$ is still positive definite and
  is therefore a metric. Moreover, $\tilde{g}$ creates the same
  gradient flow, because for every $\mu \in \dP$ and any $v\in T$
  \begin{equation*}
    \tilde{g}|_{\mu}(\mu Q,v) = \tilde{g}|_{\mu}(e_1,v)
    = g|_{\mu}(e_1,v) = g|_{\mu}(\mu Q,v).\qedhere
  \end{equation*}
\end{proof}

For the characterisation at the equilibrium distribution, we use the
assumed differentiability of $g$ and $\fF$.
\begin{proof}[Proof of \Cref{thm:pi}]
  The equality of the flow associated to the semigroup $e^{tQ}$ and
  the gradient flow implies that the time derivatives of both
  evolutions agree at every state $\pi +h$ with $h \in T$. This means
  that, for all $v \in T$,
  \begin{equation}
    \label{eq:general-cond}
    g|_{\pi+h}((\pi+h)Q,v) = - \dd \fF|_{\pi+h}(v).
  \end{equation}
  Since $\pi Q=0$, this implies that $\dd \fF$ must vanish at
  $\pi$. Moreover, it simplifies \Cref{eq:general-cond} to
  \begin{equation*}
    g|_{\pi+h}(hQ, v) = - \dd \fF|_{\pi+h}(v).
  \end{equation*}
  Let $M = \dd^2 \fF|_{\pi}$, then the Taylor
  expansion around $h=0$ shows by the assumed regularity of $g$ and
  $\fF$ that
  \begin{equation*}
    g|_{\pi}(hQ,v) + O(\|h\|^2) = -M(h,v) + O(\|h\|^2).
  \end{equation*}
  As this holds for arbitrary $h\in T$, the linear terms must
  agree. Hence,
  \begin{equation}
    \label{eq:linear}
    g|_\pi(wQ,v) = -M(w,v) \quad \forall v,w\in T.
  \end{equation}

  Furthermore, we claim that the restriction of $Q$ to $T$ defines an
  automorphism on $T$. Since $Q$ preserves the probability mass
  (i.e. $\sum_{i\in\spX} Q_{ji} = 0$ for $j\in \spX$), its range is
  inside $T$. If $Q$ was not an automorphism, a $v \in T$ satisfying
  $vQ = 0$ would exist by the Rank-Nullity Theorem. But then, for
  small enough $\alpha$, also $\pi+\alpha v$ would be a stationary
  state, contradicting the irreducibility of the Markov chain.

  Hence \Cref{eq:linear} determines the value of $g|_{\pi}$ for all
  arguments, which proves part \ref{item:g-unique} of the theorem.

  Given $g|_{\pi}$ and $M$, \Cref{eq:linear} determines $vQ\cdot w$
  for all $v,w\in T$, because $g|_{\pi}$ is a positive form. By the
  mass conservation $vQ \in T$, so that this determines $vQ$ for all
  $v \in T$. Since $\pi Q = 0$, this determines $Q$ and shows part
  \ref{item:q-dependence}.

  In order to prove part \ref{item:q-real-diag} let $f_1,\dots ,f_d$
  be a basis of $T$ which is orthonormal under $g|_\pi$,
  i.e. $g|_{\pi}(f_i,f_j) = \delta_{ij}$. Let $\bar{Q}$ be the matrix
  corresponding to the generator $Q$ in this basis, i.e. $f_i Q =
  \sum_{j=1}^{d} \bar{Q}_{ij} f_j$ for $i=1,\dots ,d$. Also let $\bar{M}$ be
  the matrix corresponding to $M$ in this basis, i.e. $\bar{M}_{ij} =
  M(f_i,f_j)$. Then \Cref{eq:linear} becomes
  \begin{equation*}
    x \bar{Q} \cdot y = - x \bar{M} \cdot y
    \quad \forall x,y \in \R^d.
  \end{equation*}
  Hence $\bar{Q} = -\bar{M}$. Since the partial derivatives of $\fF$
  commute, $\bar{M}$ is symmetric.  Therefore, $\bar{Q}$ is real
  diagonalisable and $\bar{Q}$ has $d$ real eigenvectors in
  $T$. Since $\pi$ is another eigenvector of $Q$ not in $T$, this
  implies that $Q$ is real diagonalisable, which is the statement of
  part \ref{item:q-real-diag}.

  For the converse of the theorem, we assume that $Q$ is diagonalisable and
  we need to construct a suitable functional and metric on $\dP$.  For
  this, fix eigenvectors $\pi,f_1,\dots ,f_d$ of $Q$ with eigenvalues
  $0,\lambda_1,\dots , \lambda_d$.  Since $\pi$ is the only stationary
  distribution, $\lambda_i \not = 0$ for $i=1,\dots ,d$. As $Q$ maps
  into $T$, this shows that $f_1,\dots ,f_d$ is a basis of $T$. Define
  the constant Riemannian metric $g$ on $\dP$ by
  \begin{equation*}
    g(f_i,f_j) = \delta_{ij},
  \end{equation*}
  and the functional $\fF:\dP\mapsto \R$ by
  \begin{equation*}
    \fF(\pi+\sum_{i=1}^d a_i f_i) = \frac 12 \sum_{i=1}^d (-\lambda_i) a_i^2.
  \end{equation*}
  Then at any state $\mu = \pi + \sum_{i=1}^d a_i f_i \in \dP$ we
  have, for $j=1,\dots,d$,
  \begin{equation*}
    g(\mu Q, f_j) = \lambda_j a_j \text{ and } -\dd \fF|_{\mu}(f_j) =
    \lambda_j a_j.
  \end{equation*}
  Hence the flow associated to the semigroup $e^{tQ}$ and the gradient
  flow agree, because their time derivatives agree for every
  probability distribution $\mu \in \dP$.

  For the remaining part \ref{item:q-time-reversible}, the functional
  $\fF$ is the relative entropy $\relH$.  The second derivative
  $\dd^2 \relH|_{\pi}$ of $\relH$ at $\pi$ is given by
  \begin{equation*}
    M(w,v) = \sum_{\alpha\in\spX} \frac{w_\alpha v_\alpha}{\pi_{\alpha}}
  \end{equation*}
  for $v,w \in T$.  Let $\Pi\in\R^{\spX\times\spX}$ be the diagonal
  matrix with diagonal entries $(\pi)_{i\in\spX}$. Then, by the
  calculated form of $M$, we have $v\Pi^{-1}\cdot w = M(v,w)$ for all
  $v,w \in T$.

  Over $T$, the metric $g$ has an inverse $b$ at $\pi$ which is
  defined by $g|_{\pi}(v,w b) = v \cdot w$ for $v,w \in T$ and is a
  positive definite symmetric automorphism on $T$. Define the
  symmetric matrix $a \in \R^{\spX \times \spX}$ by $v a = v b$ for
  $v \in T$ and $\mathbf{1} a = 0$, where $\mathbf{1}$ is the vector
  with all entries $1$.

  Since $b$ is an automorphism on $T$, \Cref{eq:linear} implies
  $g(vQ, ua) = {- v \Pi^{-1} \cdot (u a)}$ for all $u,v \in T$. By the
  construction of $a$, this shows that for all $u,v \in T$
  \begin{equation}
    \label{eq:reversible}
    vQ \cdot u = - v \Pi^{-1}a \cdot u.
  \end{equation}

  Since $\pi Q$ and $\pi \Pi^{-1} a$ both vanish, \Cref{eq:reversible}
  also holds for $v = \pi$. Hence it holds for all $v \in \R^{\spX}$
  which shows $Q \cdot u = -\Pi^{-1} a \cdot u$ for all $u \in T$.

  By the conservation of probability $Q\cdot \mathbf{1} = 0$ and by
  the symmetry of $a$ also $a \cdot \mathbf{1} = 0$, so that
  $Q \cdot u = -\Pi^{-1} a \cdot u$ holds for all $u \in
  \R^{\spX}$. This shows
  \begin{equation*}
    \Pi Q = -a.
  \end{equation*}
  Since $a$ is symmetric, this shows that $Q$ satisfies the detailed
  balance equation, i.e. that the Markov chain is time-reversible.
\end{proof}

The remaining \Cref{thm:relH} is reduced to the following lemma.
\begin{lemma}
  \label{lemma:relH}
  Assume the hypothesis of \Cref{thm:relH}.  Then $\fF$ and $\relH$
  have a minimum at $\pi$, and the derivatives $M=\dd^2 \fF|_{\pi}$
  and $N=\dd^2 \relH|_{\pi}$ are positive non-degenerate forms on
  $T$. For every $v\in T$ vanishing in exactly one state, there exists
  $\alpha_{v} \in \R^+$ such that $M(v,\cdot) = \alpha_{v} N(v,\cdot)$.
\end{lemma}
With this lemma the theorem can be proved.
\begin{proof}[Proof of \Cref{thm:relH}]
  By the previous lemma $\fF$ and $\relH$ have a minimum at $\pi$ so
  that $\dd \fF_{\pi} = \dd \relH_{\pi} = 0$.

  If $d=1$, then $M$ and $N$ from the lemma correspond to positive
  scalars so that there exists a positive scalar $\alpha$ satisfying
  the required relation $M=\alpha N$.

  Hence, assume $d\ge 2$. For two states $v,\tilde{v} \in T$ vanishing
  only in one common state, the constants $\alpha_v$ and
  $\alpha_{\tilde{v}}$ of the lemma must agree. If $v$ and $\tilde{v}$
  are linearly dependent, then this follows directly from the
  bilinearity of $M$ and $N$. Otherwise, for small enough
  $\lambda \in \R$, also $v+\lambda \tilde{v}$ is in $T$ and vanishing
  in exactly one state. Hence the lemma applies to this state and by
  linearity follows
  \begin{equation*}
    \alpha_{v+\lambda\tilde{v}} N(v+\lambda \tilde{v},\cdot)
    = \alpha_{v} N(v,\cdot) + \lambda \alpha_{\tilde{v}} N(\tilde{v},
    \cdot).
  \end{equation*}
  Since $N$ is non-degenerate, the functionals $N(v,\cdot)$ and
  $N(\tilde{v},\cdot)$ are linearly independent, so that the equation
  implies $\alpha_{v+\lambda{\tilde{v}}} = \alpha_{v}$ and
  $\alpha_{v+\lambda{\tilde{v}}} = \alpha_{\tilde{v}}$ and thus, as
  claimed, $\alpha_{v} = \alpha_{\tilde{v}}$.

  Hence, for $i \in \spX$, there exists $\bar{\alpha}_i$ such that
  $M(v,\cdot) = \bar{\alpha}_i N(v,\cdot)$ holds for $v\in T$ with
  $v_i = 0$ and $v_j \not = 0$ for $j \not = i$. By continuity of
  $M(v,\cdot)$ and $\bar{\alpha}_i N(v,\cdot)$ with respect to $v$, the
  result also holds for all $v \in T$ with $v_i = 0$.

  If $d \ge 3$ then for $i,j \in \spX$ there exists $v \in T$ with
  $v_i=v_j=0$ which implies that all $\bar{\alpha}_i$ for $i\in \spX$
  agree.
  If $d=2$, we can consider $v=(1\; {-1}\; 0)$ and
  $\tilde{v}=(0\; 1\; {-1})$ with $\alpha_v = \bar{\alpha}_2$ and
  $\alpha_{\tilde{v}} = \bar{\alpha}_0$. Then $v+\tilde{v}$ satisfies
  the lemma as well so that as before
  $\bar{\alpha}_0 = \bar{\alpha}_2$. Likewise,
  $\bar{\alpha}_0 = \bar{\alpha}_1$ and all $\bar{\alpha}_i$ agree.

  The common value $\alpha$ of the $\bar{\alpha}_i$ for $i \in \spX$
  is the claimed constant satisfying
  $\dd^2 \fF|_{\pi} = \alpha \dd^2 \relH|_{\pi}$. Since the
  derivatives $\dd^2 \fF|_{\pi}$ and $\dd^2 \relH|_{\pi}$ are positive
  forms, the constant $\alpha$ must be positive.
\end{proof}
The remaining lemma is proved by considering suitable Markov chains.
\begin{proof}[Proof of \Cref{lemma:relH}]
  Since $\pi$ has positive mass for every state $i \in \spX$, there
  exists a time-reversible irreducible Markov chain with stationary
  state $\pi$. Consider such a Markov chain with generator $Q$ and let
  $g$ be a metric such that the flow associated to the Markov
  semigroup $e^{tQ}$ is the gradient flow of $\fF$ under $g$. Then for
  any initial probability state $\mu \in \dP$ the value
  $\fF(\mu e^{tQ})$ is decaying as $t$ increases. Hence $\fF$ must
  have a minimum at $\pi$. As in the proof of \Cref{thm:char},
  \Cref{eq:linear} must hold and $g|_{\pi}$ and $Q$ are
  non-degenerated over $T$. Hence $M = \dd^2 \fF|_{\pi}$ must be
  non-degenerated and positive, because $\fF$ has a minimum at
  $\pi$. Since $\relH$ satisfies the assumptions imposed on $\fF$,
  this also holds for $\relH$.

  If $d=1$, there exists no such state $v \in T$. Thus assume
  $d \ge 2$ henceforth. Moreover, label the states such that $v_d = 0$
  and identify $M$, $N$ and $g|_{\pi}$ with the corresponding matrix,
  i.e.  $M(v,w) = vM \cdot w$, and $N(v,w) = vN \cdot w$, and
  $g|_{\pi}(v,w) = v g|_{\pi} \cdot w$ for all $v,w \in T$.

  Then, \Cref{eq:linear} implies $M^{-1} Q = - g|_{\pi}^{-1}$. Hence
  $M^{-1}Q$ is symmetric and thus $M^{-1} Q M = Q^T$. Likewise
  $N^{-1} Q N = Q^{T}$ and thus $M N^{-1} Q N M^{-1} = Q$.  Therefore,
  if $v$ is an eigenvector of $Q$, then $v M N^{-1}$ is again an
  eigenvector of $Q$ with the same eigenvalue.

  We finish the proof of the lemma by showing that $v$ and
  $w := v M N^{-1}$ are proportional. For this, we show $w_d = 0$ and
  $w_i v_j = v_i w_j$ for all $i,j = 0,1,\dots ,d-1$ by constructing
  suitable Markov chains whose generator $Q$ has the eigenvector $v$.

  For the first part, let $\lambda > 0$ and consider the
  time-reversible irreducible Markov chain with stationary state $\pi$
  and generator matrix
  \begin{equation*}
    Q = \Pi^{-1}
    \begin{pmatrix}
      -\pi_0 \lambda & 0 & 0 & \hdots & 0 & \pi_0 \lambda \\
      0 & -\pi_1 \lambda & 0 & \hdots & 0 & \pi_1 \lambda \\
      0 & 0 & -\pi_2 \lambda & \hdots & 0 & \pi_2 \lambda \\
      \vdots & \vdots & \vdots & \ddots & \vdots & \vdots \\
      0 & 0 & 0 & \hdots & -\pi_{d-1} \lambda & \pi_{d-1} \lambda \\
      \pi_0 \lambda & \pi_1 \lambda & \pi_3 \lambda & \hdots &
      \pi_{d-1} \lambda & -(1-\pi_d) \lambda
    \end{pmatrix}.
  \end{equation*}
  The eigenspace for $-\lambda$ is $\{u\in T : u_d=0\}$ and contains
  $v$. Hence, $w_d=0$.

  By further relabelling the states, it suffices to show for the
  second case $w_0 v_1 = v_0 w_1$.  For this, let $\lambda > 0$ and
  consider the Markov chain with generator
  \begin{align*}
    Q &= \Pi^{-1}
    \begin{pmatrix}
      -\pi_0\lambda - \beta_0 & \mu & 0 & \hdots & 0 & \pi_0\lambda + \beta_0-\mu \\
      \mu & -\pi_1\lambda -\beta_1 & 0 & \hdots & 0 & \pi_1\lambda + \beta_1-\mu \\
      0 & 0 & -\pi_2\lambda & \hdots & 0 & \pi_2 \lambda \\
      \vdots & \vdots & \vdots & \ddots & \vdots & \vdots \\
      0 & 0 & 0 & \hdots & -\pi_{d-1}\lambda & \pi_{d-1} \lambda \\
      \pi_0\lambda+\beta_0-\mu & \pi_1\lambda+\beta_1-\mu & \pi_2 \lambda & \hdots &
      \pi_{d-1} \lambda & \gamma
    \end{pmatrix}
  \end{align*}
  where
  \begin{align*}
    \beta_0 = \mu \frac{\pi_0 v_1}{v_0 \pi_1},\qquad
    \beta_1 = \mu \frac{\pi_1 v_0}{v_1 \pi_0},\qquad
    \gamma = -(1-\pi_d) \lambda - \mu \frac{\pi_0 v_1}{v_0 \pi_1}
    -\mu \frac{\pi_1 v_0}{v_1 \pi_0} + 2\mu.
  \end{align*}
  By choosing $\mu>0$ small enough, this defines an irreducible
  time-reversible Markov chain with stationary state $\pi$ and
  $vQ = -\lambda v$. Hence, $w$ is again an eigenvector with
  eigenvalue $-\lambda$. From the first component, we find the
  required ratio
  \begin{equation*}
    -\lambda w_0 = -\lambda w_0 - \mu \frac{v_1 w_0}{v_0\pi_1} + \mu
    \frac{w_1}{\pi_1}
    \Rightarrow
    v_1w_0 = v_0 w_1. \qedhere
  \end{equation*}
\end{proof}





\ACKNO{The author would like to thank James Norris for many helpful
discussions and suggestions.}


\end{document}